\documentclass[a4paper,11pt]{article}

\usepackage{amsfonts} 
\usepackage{amssymb} 
\usepackage{amsmath,amsthm} 

\newtheorem{theoreme}{Theorem}[section] 
\newtheorem{lemma}[theoreme]{Lemma} 
\newtheorem{prop}[theoreme]{Proposition} 
\newtheorem{cor}[theoreme]{Corollary} 
\newtheorem{definition}[theoreme]{Definition}

\def \cqfd {\ { \hspace{\stretch{9} } $\square$}} 
\def \Rk {\ {\bf Remark.} } 
\def \Rks {\ {\bf Remarks.}} 
\def \sm {\setminus } 
\newcommand{\be}{\begin{enumerate}}  \newcommand{\ee}{\end{enumerate}} 
\newcommand{\bi}{\begin{itemize}}  \newcommand{\ei}{\end{itemize}} 
\newcommand{\bd}{\begin{description}}  \newcommand{\ed}{\end{description}} 
\newcommand{\ra}{\rightarrow} 
\newcommand{\lan}{\langle}   \newcommand{\ran}{\rangle}

\newcommand{\comment}[1]{}

\def \R {\mathbb{R}}

\def \harm {\mathcal{H}}

\numberwithin{equation}{section}  
\renewcommand{\phi}{\varphi} 
\renewcommand{\epsilon}{\varepsilon} 
 
\title{$L^2$-cohomology of negatively curved K\"ahler manifolds of finite volume}
\author{N. Yeganefar}
\date{February 2, 2004} 

\begin{document} 
\maketitle

\begin{abstract}
We compute the space of $L^2$ harmonic forms (outside the middle degrees) on negatively curved K\"ahler manifolds of finite volume.
\end{abstract}

\section{Introduction}
Let $(M,g)$ be a Riemannian manifold. We denote by $\harm ^k(M)$ the space of $L^2$ harmonic $k-$forms on $M$. When $M$ is compact, the Hodge-de Rham theorem states that this space is finite dimensional and in fact isomorphic to the k$^{th}$ real cohomology group of $M$. When $M$ is not compact, there is no general theorem which identifies $\harm ^k(M)$ with a topologically defined group. However, for a number of interesting non compact manifolds, Hodge type theorems have already been proved, see e.g. the work of Atiyah-Patodi-Singer \cite{APS} on manifolds with cylindrical ends, of S. Zucker \cite{Z} on locally symmetric manifolds, of R. Mazzeo and Mazzeo-Phillips \cite{M,MP} on (asymptotically) hyperbolic manifolds, of G. Carron \cite{C2} on manifolds with flat ends, etc. 

There are also  cases for which we know that $\harm ^k(M)$ is finite dimensional, but no topological interpretation is available. For example, according to J. Lott \cite{L2}, if $M$ is a complete manifold of finite volume and pinched negative curvature, then its spaces of $L^2$ harmonic forms are all finite dimensional. Lott also raised the problem of computing the dimensions of these spaces. As a first step towards solving this problem, we proved the following result in \cite{Y} (the constant curvature case was already obtained by Zucker \cite{Z} and Mazzeo-Phillips \cite{MP}): 
\begin{theoreme}[\cite{Y}]\label{cohomology}
Let $(M^n,g)$ be a complete $n-$dimensional Riemannian manifold of finite volume and with pinched negative sectional curvature $K$: there exist two constants $0<a<b$ such that $-b^2\leq K\leq -a^2$. Assume that $na-(n-2)b>0$. Then we have the isomorphisms
$$\mathcal{H}^k(M)\simeq \left\{ \begin{array}{lll}
     H^k(M), & \textrm{if $k<(n-1)/2$,}\\
{\rm Im} (H^{n/2}_c(M)\ra H^{n/2}(M)), & \textrm{if $k=n/2$,}\\
     H^k_c(M), & \textrm{if $k> (n+1)/2.$}
     \end{array} \right. $$
If moreover the curvature is constant ($a=b$) and the dimension $n$ is odd, then we have  $\mathcal{H}^{(n\pm 1)/2}(M)\simeq {\rm Im} (H^{(n\pm 1)/2}_c(M)\ra H^{(n\pm 1)/2}(M))$. 
\end{theoreme}
This means that if the curvature is sufficiently pinched, then we have a topological interpretation of $\mathcal{H}^*(M)$. We also gave examples showing that our theorem is sharp with respect to the pinching constants: as soon as $a$, $b$ and $n$ satisfy $na-(n-2)b\leq 0$, there is a manifold with the correct geometry but for which the conclusion of theorem \ref{cohomology} doesn't hold; similarly, for each $a<b$ and for each $n$ odd, there is a $n$-dimensional manifold with finite volume and curvature pinched between $-b^2$ and $-a^2$ for which the topological interpretation in degrees $(n\pm 1)/2$ fails. However, it follows from Zucker's work \cite{Z} that for finite volume quotients of the complex hyperbolic space the conclusion of theorem \ref{cohomology} holds, although these manifolds do not satisfy the pinching assumption of this theorem. The natural problem is then to try to find an analog of theorem \ref{cohomology} for K\"ahler manifolds. This is also a problem raised by S.T. Yau \cite[Question 61]{Yau}. Our main result in this paper partially solves this problem:
\begin{theoreme}\label{kaehler}
 Let $(M^{2m},g)$ be a complete  K\"ahler manifold of real dimension $n=2m$, with finite volume and  pinched negative sectional curvature $-b^2\leq K\leq -a^2<0$. Then we have the isomorphisms
$$\mathcal{H}^k(M)\simeq \left\{ \begin{array}{ll}
     H^k(M), & \textrm{if $k<m-1$,}\\
     H^k_c(M), & \textrm{if $k>m+1.$}
     \end{array} \right. $$
Moreover, there is an injection of $\mathcal{H}^{m-1}(M)$ in $H^{m-1}(M)$ and a surjection of $H^{m+1}_c(M)$ onto $\mathcal{H}^{m+1}(M)$.
\end{theoreme}
We believe that there are isomorphisms $\mathcal{H}^{m-1}(M)\simeq H^{m-1}(M)$, $\mathcal{H}^{m+1}(M)\simeq H_c^{m+1}(M)$, and $\mathcal{H}^{m}(M)\simeq {\rm Im} (H^{m}_c(M)\ra H^{m}(M))$, but we have not been able to show this.

Let us say a few words about the proof of Theorem \ref{kaehler}. Our strategy is similar to ideas developed in \cite{Y,Y2} (and this is also implicit in the paper of Donnelly-Fefferman \cite{DF}). Thus, we proceed as follows:
\be
\item[i)] We show that zero is not in the essential spectrum of the Laplacian. \item[ii)] As is well-known, this implies that there is an exact sequence which relates the relative cohomology of any bounded open subset $D\subset M$ to the $L^2$-cohomology of $M$ and the absolute $L^2$-cohomology of $M\sm D$.
\item[iii)] We show that the absolute $L^2$-cohomology of $M\sm D$ vanishes (for suitable degrees) and conclude \textit{via} the exact sequence.  
\ee  
More specifically, using Busemann functions, we first prove that outside a compact subset, the K\"ahler form of the metric is $d$(bounded) (in the terminology of Gromov \cite{Gro}); this implies i). The difficult step is then to get iii). For this, our main tool is the construction of a new complete K\"ahler metric on each end of $M$, by using again Busemann functions.

The paper is organized as follows. We begin with a review of some classical facts about the space of $L^2$ harmonic forms and its links with $L^2$-cohomology (see Section 2). In section 3, we introduce Gromov's notion of K\"ahler hyperbolicity, and say a few words about the geometry of negatively curved manifolds of finite volume. We then combine this material to study the essential spectrum of negatively curved K\"ahler manifolds of finite volume. Section 4 is more technical: we construct the new K\"ahler metric which was mentioned above, and we study its main properties. Finally, we prove Theorem \ref{kaehler} in Section 5.\\

\textbf{Acknowledgements.} This paper is taken form my PhD thesis. Just one word for my PhD supervisor G. Carron: merci pour tout, Gilles! I also thank J. Br\"uning and his team, especially G. Marinescu, for their hospitality at the Humboldt Universit\"at (Berlin), while part of this work was written.

\section{$L^2$ harmonic forms and $L^2$-cohomology}
We present in this section classical results about $L^2$ harmonic forms. For more details, the reader may consult the works of J. Lott \cite{L1} or of G.  Carron \cite{C1,C2}.
\subsection{$L^2$ harmonic forms}
Let $M$ be a manifold, and let $d$ denote the differential which acts for example on smooth compactly supported $k-$forms:
  $$d :  C^\infty _0(\Lambda ^kT^*M)\to C^\infty _0(\Lambda ^{k+1}T^*M).$$
If $M$ is endowed with a Riemannian metric $g$, we can consider the formal adjoint of $d$, to be denoted by $\delta$:
  $$\delta : C^\infty _0(\Lambda ^{k+1}T^*M)\to C^\infty _0(\Lambda ^kT^*M)$$ is defined by
  $$\forall \alpha , \beta \in C^\infty _0, \lan \delta \alpha ,\beta \ran _{L^2}=\lan \alpha ,d\beta \ran _{L^2}.$$
\begin{definition}
The space $\mathcal{H}^k(M)$ of $L^2$ harmonic  $k-$forms of $(M,g)$ is
$$\mathcal{H}^k(M)=\{\alpha \in L^2(\Lambda ^kT^*M)/\, d\alpha =\delta \alpha =0\},$$
where $L^2(\Lambda ^kT^*M)$ denotes the space of square integrable $k$-forms, and the equations are first understood in the weak sense.
\end{definition}
Denote by $\Delta =(d+\delta)^2$ the Laplace operator and by $\Delta _k$ its restriction acting on $k-$forms. If the metric is complete, then we can integrate by parts and get 
$$\mathcal{H}^k(M)=\rm{Ker}_{L^2}(\Delta _k).$$

\subsection{$L^2$-cohomology}
 \subsubsection{Complete manifolds}
We now first assume that $M$ is complete.
\begin{definition}\label{L2coh}
The $k^{th}$ (reduced) $L^2$-cohomology space of $M$ is
 $$H^k_2(M)=\{\alpha \in L^2(\Lambda ^kT^*M)/\, d\alpha =0\}/\overline{dC^\infty _0(\Lambda ^{k-1}T^*M)}^{L^2} ,$$ 
where $\overline{dC^\infty _0(\Lambda ^{k-1}T^*M)}^{L^2}$ means that we take the closure of \\$dC^\infty _0(\Lambda ^{k-1}T^*M)$ in $L^2$.
\end{definition}
We have the following Hodge-de Rham-Kodaira decomposition \cite[th\'eor\`eme 24]{dR}:
$$L^2(\Lambda ^kT^*M)= \mathcal{H} ^k(M) \oplus
  \overline{dC^\infty _0(\Lambda ^{k-1} T^*M)} \oplus
  \overline{\delta C^\infty _0(\Lambda ^{k+1} T^*M)},$$
and moreover
  $$\{\alpha \in L^2(\Lambda ^kT^*M)/\, d\alpha =0\}=\mathcal{H} ^{k}(M) \oplus
  \overline{dC^{\infty}_{0}(\Lambda ^{k-1} T^*M)}.$$
It follows that $$H^k_2(M)\simeq \mathcal{H}^k(M),$$
and from now on we won't really distinguish these spaces.

A closely related space is the unreduced $L^2$-cohomology space, which is by definition the cohomology of the $L^2$ de Rham complex (without taking the closure of the image of $d$). Reduced and unreduced $L^2$-cohomology coincide if and only if the $L^2$ image of $d$ is closed. This is the case if for example zero is not in the essential spectrum. For us, $L^2$-cohomology means reduced $L^2$-cohomology. Note however that by Corollary \ref{essentiel} below, we consider in this paper situations where reduced and unreduced $L^2$-cohomology are the same.

  \subsubsection{Manifolds with boundary}
We assume here that $(M,g)$ is a manifold with smooth compact boundary, and that it is metrically complete. We can then define absolute (or relative) $L^2$-cohomology groups.

Denote by $C^\infty _b(\Lambda ^kT^*M)$ the space of smooth $k-$forms on $M$ whose support is bounded (this support can meet the boundary for elements of $C^\infty _b$, but not for elements of $C^\infty _0$). The $k^{th}$ absolute $L^2$ cohomology group is defined by:
$$H^k_2(M)=(\delta C^\infty _0(\Lambda ^{k+1}T^*M))^\perp/\overline{dC^\infty _b(\Lambda ^{k-1}T^*M)} .$$
We have the following identification of $H^k_2(M)$ with a space of harmonic forms satisfying an  absolute boundary condition on $\partial M$:
$$H^k_2(M)\simeq \{\alpha \in L^2(\Lambda ^kT^*M)/\, d\alpha =\delta \alpha =0,\, i_\nu\alpha =0\},$$
where $\nu$ is a normal vector field on the boundary, and $i_\nu\alpha$ denotes interior product.

\subsection{$L^2$-cohomology and essential spectrum}\label{spectre}
Recall that the discrete spectrum of the Laplacian $\Delta$ is by definition the set of eigenvalues which are of finite multiplicity and isolated in the spectrum; the essential spectrum is the complement in the spectrum of the discrete spectrum. It is well known (see \cite{G,D,A,Ba}) that the essential spectrum does only depend on the geometry at infinity. More precisely, we have the following characterization: zero is not in the essential spectrum of $\Delta _k$ if and only if we have a Poincar\'e inequality at infinity, i.e. there exist a compact subset $D\subset M$ and a positive constant $C>0$ such that
$$ \forall \alpha \in C^{\infty}_{0}(\Lambda ^{k}
T^{*}(M\sm \overline{D})),\, C\Vert \alpha \Vert
_{L^{2}}\leq \Vert \Delta _k\alpha \Vert_{L^{2}}.$$ 

By a result due to J. Lott \cite{L1}, the finiteness of the dimension of $\mathcal{H} ^k(M)$ depends also only on the geometry at infinity. It is therefore natural to look for links between the $L^2$-cohomology of a manifold, the toplogy of this manifold and its geometry at infinity. The following proposition provides such links for the case when zero is not in the essential spectrum of the Laplacian $\Delta _k$. 
\begin{prop}\label{se}
Let $(M,g)$ be a complete Riemannian manifold and let $D$ be a bounded open subset in $M$ with regular boundary. Assume that for some $k$, zero is not in the essential spectrum of $\Delta _k.$ Then we have the following exact sequence
$$H^{k-1}_c(D)\buildrel e\over \to H^{k-1}_2(M)\buildrel r\over \to H^{k-1}_2(M\setminus D)
\buildrel b \over \rightarrow H^k_c(D)\buildrel e \over \rightarrow H^k_2(M)$$
$$\buildrel r \over \rightarrow H^k_2(M\sm D)  \buildrel b \over \rightarrow H^{k+1}_c(D) \buildrel e \over \rightarrow H^{k+1}_2(M) \buildrel r \over \rightarrow H^{k+1}_2(M\sm D) \buildrel b \over \rightarrow H^{k+2}_c(D),$$ 
where $e$ is extension by zero, $r$ is restriction, and $b$ is the coboundary operator.
\end{prop}
This is a well-known result. For the definition of the maps $r$, $e$, $b$, and for the proof of this proposition, the reader is referred to \cite{Y} (see also \cite{C1,C2}).

\section{Manifolds of negative curvature and Gromov's K\"ahler hyperbolicity}

 \subsection{K\"ahler hyperbolicity}
In \cite{Gro}, M. Gromov introduced a new notion of hyperbolicity. More specifically, following  Gromov, we say that a differential form $\alpha$ on a Riemannian manifold is "$d$(bounded)" if there exists a bounded form $\beta$ such that $\alpha =d\beta$. We say that $\alpha$ is $\tilde{d}$(bounded) if its lift to the universal covering is $d$(bounded). A compact complex manifold is then called K\"ahler hyperbolic if it admits a K\"ahler metric whose fundamental form is $\tilde{d}$(bounded). The most important examples are compact complex manifolds which admit a K\"ahler metric of negative curvature. 
Gromov used this notion to show that on the universal covering of a K\"ahler hyperbolic manifold, there are no (non trivial) $L^2$ harmonic forms outside the middle degree (he also showed that the space of $L^2$ harmonic forms in the middle degree is infinite dimensional). A key argument in Gromov's proof is the following result 
\begin{theoreme}[Gromov]\label{Gromov}
Let $(M^{2m},\omega )$ be a complete K\"ahler manifold of real dimension $2m$. Assume that the fundamental form $\omega$ is $d$(bounded), i.e. there exists a bounded $1-$form $\theta$ such that $\omega =d\theta$. Then if $p\neq m$ and if $\alpha$ is a $p-$form in the domain of the Laplacian, we have the estimate 
$$c_m \Vert \theta \Vert _{L^\infty }^{-2} \Vert \alpha \Vert _{L^2} \leq \Vert \Delta \alpha \Vert _{L^2} ,$$ where $c_m>0$ is a constant which depends only on $m$.
In particular, we have $\mathcal{H}^p(M)=0$ for $p\neq m$.
\end{theoreme}
\Rk If we suppose that we have $\omega =d\theta$ outside a compact subset, then Gromov's argument can be applied and shows that we have a Poincar\'e inequality at infinity for $\Delta _p$, $p\neq m$ (see the proof of \cite[Theorem 1.4.A]{Gro}). Thus zero is not in the essential spectrum of $\Delta _p$, $p\neq m$.

\subsection{Manifolds with pinched negative curvature and finite volume}
In this section, $(M,g)$ is complete manifold of finite volume and with pinched negative sectional curvature $K$, i.e. there are constants $0<a<b$ such that $$-b^2\leq K\leq -a^2.$$
In the next subsection, we will use Gromov's results to study spectral properties of these manifolds. But let us recall for the moment some standard facts about the topology and geometry of these manifolds (see \cite{E} and \cite{HI}). First, $M$ has a finite number of ends, and one has $M=M_0\cup E_i$, where $M_0$ is a compact manifold with boundary, and the $\partial E_i$'s are the components of $\partial M_0$. To each ray of $E_i$, we can associate a Busemann function $r_i$ which is
 a priori only $C^2$-smooth. Two such functions are equal up to an additive
 constant. $E_i$ is $C^2-$diffeomorphic to $\R^+\times \partial E_i$. Moreover, the slices $\{ t\}\times \partial E_i$ are the level sets of a Busemann function. Finally, the metric on each end $E_i$ has the following form
$$g=dr_i^2+h_{r_i},$$
where $h_{r_i}$ is a family of metrics on the compact manifold $\partial E_i$, and satisfies  $e^{-r_i}h_0\leq h_{r_i}\leq e^{-ar_i}h_0.$

Assume now that the metric $g$ is moreover a K\"ahler metric. Then it was observed by Greene-Wu \cite{GrW} and Siu-Yau \cite{SY} that each function $-r_i$ is strictly plurisubharmonic on the corresponding end. More precisely, on each end $E_i$, we have
$$ag\leq -\sqrt{-1}\partial \overline{\partial}r_i\leq bg.$$
 \subsection{Essential spectrum}
In this section, $M$ is a K\"ahler manifold of finite volume and with pinched negative curvature. On such a manifold, the K\"ahler form $\omega$ cannot be $d$(bounded). Namely, on any K\"ahler manifold the fundamental form $\omega$ has constant length and is harmonic (because it is parallel). As the volume is finite, $\omega$ is a non trivial $L^2$ harmonic form. But because of the Hodge-de Rham decomposition, $\omega$ cannot be written as the differential of a bounded (hence $L^2$) form. However, we have the following crucial lemma, which is probably well known:
\begin{lemma}\label{hyp}
Let $(M,g)$ be a complete K\"ahler manifold of finite volume and with pinched negative curvature $-b^2\leq K\leq -a^2<0$. Then outside a compact subset, its K\"ahler form is $d$(bounded). More specifically, there exist a bounded open subset $D\subset M$ and a bounded and continuous $1-$form $\theta$ such that we have $\omega =d\theta$ in the weak sense on $M\sm D$.
\end{lemma}
\begin{proof}
Without loss of generality, we may assume that $M$ has only one end $E$. Let $r$ be a Busemann function associated to $E$ such that $E$ is diffeomorphic to a product $(0,\infty)\times \Sigma$, with $\Sigma$ a closed connected manifold. Let $\phi _t$ be the flow of the gradient $\nabla r$ of $r$. If $x=(s,y)$ is a point of $(0,\infty )\times \Sigma$, we have 
$$\phi _t(s,y )=(s+t,y ).$$ 
The properties of the differential of this flow are well understood. First, for $x\in E$, we have $d_x\phi _t(\nabla r(x))=\nabla r(\phi _t(x))$. Let $u$ be a tangent unitary vector at $x$, and orthogonal to $\nabla r$. Then $t\to d_x\phi _t(u)$ is a stable Jacobi field along the geodesic $t\to(s+t,y)$, orthogonal to $\nabla r$, and is equal to $u$ at $t=0$. As the sectional curvature is less than or equal to $-a^2$, the classical estimates on Jacobi fields \cite{HI} show that 
$$\vert d_x\phi _t(u)\vert \leq e^{-at}.$$ 
Now, if $\omega$ is the K\"ahler form of the metric, it is closed, so that by Cartan formula, we have 
$$\phi _t^*\omega -\omega =d\int _0^t \phi _s^*(i_{\nabla r}\omega )\, ds.$$ 
By using the estimate on the differential of $\phi _t$, we can take the limit as $t$ goes to infinity and get that $\omega =d\theta$ on the end $E$, with 
$$\theta =-\int _0^\infty \phi _s^*(i_{\nabla r}\omega )\, ds.$$ 
The same estimate implies that $\theta$ is bounded.
\end{proof}
Borel and Casselman \cite{BoC} studied the essential spectrum of locally symmetric manifolds of finite volume. Our next corollary is complementary to their results: 
\begin{cor}\label{essentiel}
Let $(M,g)$ be a complete K\"ahler manifold of finite volume and with pinched negative curvature. Then zero is not in the essential spectrum of the Laplacian $\Delta$ acting on $L^2$ differential forms.
\end{cor}
\begin{proof} It follows from Lemma \ref{hyp} and the remark after Gromov's theorem \ref{Gromov} that zero is not in the essential spectrum of $\Delta _p$, for $p\neq m$ (where $2m$ is the real dimension of $M$). Moreover, a classical argument shows that there is a spectral gap at zero for $\Delta _m$. It is therefore enough to show that $\harm ^m(M)$ is finite dimensional, and this is proved in \cite{L2} by Lott. 
\end{proof}
\Rk Ballmann and Br\"uning \cite{BB} (see also Donnelly and Xavier \cite{DX}) proved that if $M^{2m}$ is a $2m$-dimensional manifold of finite volume and with pinched negative curvature, but not necessarily K\"ahlerian, then its Gauss-Bonnet operator is Fredholm, provided that $ma-(m-1)b>0$ (actually, they  studied more general Dirac type operators). Moreover, they showed by exhibiting explicit examples that their pinching assumption is sharp. It follows by our corollary that in the K\"ahler case, no condition on the pinching is needed to obtain this Fredholmness property. 
\section{An auxiliary metric on the ends}
In this section, we derive some technical results which will be used in the proof of Theorem \ref{kaehler}. Thus, $(M,g)$ is a complete K\"ahler manifold of finite volume and with pinched negative curvature. We assume that $M$ has only one end $E$. Let $r$ be a Busemann function associated to $E$ such that $E$ is diffeomorphic to a product $(0,\infty)\times \Sigma$, with $\Sigma$ a closed connected manifold. Set
  $$D_0=\{ r<0\} \,\, \mathrm{and}\, \,  D_1=\{ r<1\} .$$
We construct now a complete metric on $M\sm D_0$ (sending the boundary to infinity) which has "nice" properties:
\begin{lemma}\label{metric1}
There exists on $M\sm D_0$ a complete K\"ahler metric $\tilde{g}$ such that
\be
 \item $\tilde{g}$ coincides with the initial metric $g$ on the region $\{ r\geq 2\}$,
 \item the K\"ahler form $\tilde{\omega}$ associated to $\tilde{g}$ is $d$(bounded).
\ee
\end{lemma} 
\begin{proof}
Consider a smooth function $f:(0,\infty )\to \R$ such that
$$f(t)=-\log (t) \,\textrm{on $(0,1]$ and $f'(t)=0$ on $[2,\infty )$},$$
$$f'\leq 0 \,\, \textrm{and $f''\geq 0$}.$$
Define the function $F :M\sm D_0\to \R$ by
$$F(x)=f(r(x)).$$
Recall that the Busemann function $r$ is \textit{a priori} only $C^2$, but we assume first that it is smooth. We compute $$\partial \overline{\partial}F =f'(r)\partial \overline{\partial}r+f''(r) \partial r\wedge \overline{\partial} r.$$
As $-r$ is strictly plurisubharmonic on $M\sm D_0$, the assumptions on $f$ imply that $i\partial \overline{\partial}F$ is nonnegative on $M\sm D_0$. We have a bit more: $f(t)=-\log (t)$ if $t\leq 1$, so that on  $\{0 <r\leq 1\}$ we get 
\begin{equation}\label{MD}
\partial \overline{\partial}F =-\partial \overline{\partial}(\log(r))=-\frac{\partial \overline{\partial}r}{r}+ \frac{\partial r\wedge \overline{\partial} r}{r^2},
\end{equation}
hence $i\partial \overline{\partial}F$ is positive definite and complete on $\{0 <r\leq 1\}$.  We then consider the complete K\"ahler metric $\tilde{g}$ on $M\sm D_0$ whose fundamental form $\tilde{\omega}$ is given by  
$$\tilde{\omega}=i\partial \overline{\partial}F +\omega =d(-i\partial F +\theta),$$
where we write $\omega =d\theta$ for the metric $g$, according to Lemma \ref{hyp}. We observe that $g$ and $\tilde{g}$ coincide on $\{ r\geq 2\}$ because $f$ is constant on $[2,\infty )$. 

It remains to show that $\tilde{\omega}$ is $d$(bounded), and it is enough to check that the $1-$form $-i\partial F +\theta$ is bounded with respect to $\tilde{g}$. On $\{ r>2\}$, this is clear because $\tilde{\omega }=\omega$. We now focus on the region $\{ 0<r\leq 1\}$, where $\tilde{g}$ is quasi-isometric to $i\partial \overline{\partial}F$. As $-\partial \overline{\partial}r$ is positive, equation \ref{MD} implies that
$$\vert \partial F \vert _{i\partial \overline{\partial}F} =\vert \frac{\partial r}{r} \vert _{i\partial \overline{\partial}F}\leq 1.$$
But the (pointwise) norm of $\theta$ with respect to $i\partial \overline{\partial}F$ goes to zero when $r$ is near zero, so that $-i\partial F +\theta$ is bounded with respect to $\partial \overline{\partial}F$ , and hence with respect to $\tilde{g}$.

We finish by dealing with the regularity question. We assumed that $r$ is smooth, but this is not a serious restriction. Namely, to construct our metric $\tilde{g}$, we used $r$ only on the compact subset $\{ 0\leq r\leq 2\}$. Therefore, if $r$ is not smooth, we can approximate it by a smooth function on this subset and use this approximation in the construction.
\end{proof}
We will need to study the asymptotic behaviour of the metric $\tilde{g}$ as $r$ goes to zero. We introduce first some notations. For $t\geq 0$, the restriction of the tangent bundle $TM$ to the submanifold $\{ r=t\}$ will be denoted by $TM \vert _{\{r=t \}}$. As $\{ r\geq 0\}$ is diffeomorphic to a product $[0,\infty )\times \Sigma$, we have $TM \vert _{\{r=t \}}=\R \partial/\partial r\oplus T\Sigma$. Let $J$ be the complex structure, and let $J_t$ be the restriction of $J$ to $TM \vert _{\{r=t \}}$. Set $e_1=\partial /\partial r$ and $e_2=J_0\partial /\partial r$. We choose vector fields $e_3,\ldots ,e_{2m}$ such that $e_1,e_2,e_3,\ldots ,e_{2m}$ is a local $g-$orthonormal frame of $TM\vert _{\{ r=0\} }$: this means that the $e_p$'s are local orthonormal sections of $TM\vert _{\{ r=0\} }$. Note that $e_p, p\geq 2$ are in the kernel of the $1$-form $dr$ and hence they give a local frame of $T\Sigma$, i.e. local sections which form a basis at each point. Thus, we can consider that for $t\geq 0$, the $e_p$'s form a local frame of $TM \vert _{\{r=t \}}=\R \partial/\partial r\oplus T\Sigma$ (which will not be $g-$orthonormal anymore in general). Finally, denote by $h_0$ the matrix $(-i\partial \overline{\partial}r(e_p, J_re_q))\vert _{r=0},p,q\geq 3$. With these preliminaries in mind, we can state our second technical lemma
\begin{lemma}\label{metric2}
As $t$ goes to zero, the metric $\tilde{g}$ of Lemma \ref{metric1} is quasi-isometric to the metric whose restriction to  $TM \vert _{\{r=t \}}, t>0$ is given in the frame $e_1,\ldots ,e_{2m}$ by the matrix
 $$\left( \begin{array}{cc}
  \begin{array}{cc} 1/t^2 & 0\\ 0 & 1/t^2 \end{array} & 0 \\ 0 & h_0/t \end{array} \right) .$$
\end{lemma}
\begin{proof}
It is enough to prove the lemma for the metric associated to $-i\partial \overline{\partial}\log (r)$ instead of $\tilde{g}$, because these are quasi-isometric near $r=0$ (see the proof of Lemma \ref{metric1}). For $t>0$, we have 
$$-i\partial \overline{\partial}\log (r)(e_p,J_te_q)=-i\frac{\partial \overline{\partial}r(e_p,J_te_q)}{t}+i\frac{\partial r\wedge \overline{\partial }r(e_p,J_te_q)}{t^2}.$$
We also have $$2\partial r(e_1)=dr(\partial /\partial r)-idr(J_t\partial /\partial r)=1,$$
$$2\partial r(J_te_1)=dr(J_t\partial /\partial r)+idr(\partial /\partial r)=i.$$
Hence $$-i\partial \overline{\partial}\log (r)(te_1,tJ_te_1)=-it^{-1} \partial \overline{\partial }r(te_1,tJ_te_1)+1/2=1/2+O(t).$$
As  $J_tJ_0\partial /\partial r=-\partial /\partial r+O(t)$, it follows that    $$2\partial r(e_2)=dr(J_0\partial /\partial r)-idr(J_tJ_0\partial /\partial r)=i+O(t),$$
$$2\partial r(J_te_2)=dr(J_tJ_0\partial /\partial r)+idr(J_0\partial /\partial r)=-1+O(t).$$
Finally, for $p\geq 3$, we have $dr(e_p)=0$ and $dr(J_te_p)=O(t)$, which implies that
$$\partial r(e_p)=\partial r(J_te_p)=O(t).$$ 
Consequently, in the local frame $te_1, te_2, \sqrt{t}e_3,\ldots ,\sqrt{t}e_{2m}$ ($t>0$) of $TM\vert _{\{ r=t\} }$, the matrix of $-i\partial \overline{\partial}\log (r)(.,J.)$ is
$$\left( \begin{array}{cc}
  \begin{array}{cc} 1/2+O(t) & O(t)\\ O(t) & 1/2+O(t) \end{array} & O(\sqrt{t}) \\ O(\sqrt{t}) & -i\partial \overline{\partial}r(e_p, J_te_q)+O(t) \end{array} \right) .$$
Therefore, the metric associated to $-i\partial \overline{\partial}\log (r)$ is quasi-isometric to the metric whose restriction to $TM \vert{\{ r=t\} }$ in the frame $te_1, te_2, \sqrt{t}e_3,\ldots ,\sqrt{t}e_{2m}$ is given by the matrix
$$\left( \begin{array}{cc}
  \begin{array}{cc} 1 & 0\\ 0 & 1 \end{array} & 0 \\ 0 & h_0 \end{array} \right) .$$
The conclusion of the lemma follows.
\end{proof}

\section{Proof of Theorem \ref{kaehler} and final comments}
By duality, it suffices to consider the degrees $k\geq m+1$. Without loss of generality, we assume that $M$ has only one end and we keep the notations of the previous section. By corollary \ref{essentiel}, zero is not in the essential spectrum of the Laplacian. Hence, the sequence of Proposition \ref{se} is exact for all $k$ and for the bounded subset $D_1$:
$$H^{k-1}_2(M\sm D_1)\buildrel b\over \to H^k_c(D)\buildrel e\over \to H^k_2(M)\buildrel r\over \to H^k_2(M\sm D_1).$$
To prove the theorem, it is therefore enough to show that for $j\geq m+1$, we have $H^j_2(M\sm D_1,g)=0$. The notation $H^j_2(M\sm D_1,g)$ means that we indicate the explicit dependence of the $L^2$-cohomology spaces with respect to the metric $g$. 

Let $[\alpha ]$ be a class in $H^j_2(M\sm D_1, g)$, with representative the smooth form $\alpha$. Near the boundary $\{ 1\} \times \Sigma$ of $M\sm D_1$, we write $$\alpha =\beta (r)+dr \wedge \gamma (r),$$
where $\beta$ and $\gamma$ are viewed as forms on $\Sigma$ which depend on the parameter $r$. Setting $$\phi (r)=\int _1^r \gamma (t) dt,$$ it is easily seen that $\alpha -d\phi$ doesn't have any $dr$-component. Moreover, si $d^\Sigma$ is the differential on $\Sigma$, we have
\begin{eqnarray*}
\frac{\partial }{\partial r}(\alpha -d\phi ) &=& \frac{\partial \beta }{\partial r}-d^\Sigma\gamma \\
 &=& d\alpha =0.
\end{eqnarray*}
Hence, we may replace $\alpha$ by $\alpha -d(\rho\phi)$, where $\rho$ is a smooth compactly supported function with $\rho =1$ near $\{ 1\} \times \Sigma$, and we may assume that near $\{ 1\} \times \Sigma$, the form $\alpha $ does not depend on $r$ and is without $dr$-component. We can therefore extend $\alpha$ on $M\sm D_0$ to get a form $\tilde{\alpha}$ which satisfies $\tilde{\alpha} (r)= \alpha (1)$ if $r\leq 1$. Note that $\tilde{\alpha}$ is still a closed form. 

The next observation is that $\tilde{\alpha}$ is square integrable on $M\sm D_0$ with respect to the metric $\tilde{g}$ of Lemma \ref{metric1}. Namely, by Lemma \ref{metric2}, the volume form of $\tilde{g}$ behaves like $t^{-(m+1)}$ when $t$ approches $0$, whereas the assumptions on $\tilde{\alpha}$ imply that its pointwise norm $\vert \tilde{\alpha}\vert ^2$ behaves like $O(t^{j})$. As $j\geq m+1$, $\tilde{\alpha}$ is square integrable around $t=0$.

To sum up, $\tilde{\alpha}$ is a closed $L^2$ $j$-form on $(M\sm D_0 ,\tilde{g})$, with $j\geq m+1$. But the metric $\tilde{g}$ is K\"ahler hyperbolic (see Lemma \ref{metric1}), so that  $H^j_2(M\sm D_0 ,\tilde{g})=0$ by Gromov's theorem \ref{Gromov}. It follows that there exists a form $\beta$ which is $L^2$ on $(M\sm D_0 ,\tilde{g})$ and such that $\tilde{\alpha}=d\beta$. Then we have $\alpha =d(\beta \vert _{M\sm D_1})$, where $\beta \vert _{M\sm D_1}$ is square integrable on $(M\sm D_1,g)$ (because $\tilde{g}$ and $g$ coincide if $r$ is large enough). This shows that $[\alpha ]$ is zero and that $H^j_2(M\sm D_1,g)=0$ for $j\geq m+1$. The proof of Theorem \ref{kaehler} is complete.\cqfd \\

\Rks 
\be
 \item As pointed out in the introduction, we probably have the expected conclusion for the degrees $k=m,m\pm 1$. As in \cite{Y}, this would follow if we could for example prove the vanishing of $H^m(M\sm D_1)$.
 \item We can perturb the topology and the geometry of our manifold on a compact subset, and the conclusion of Theorem \ref{kaehler} still holds (see the proof).
 \item The method of proof of Theorem \ref{kaehler} can be applied to a more general setting. Namely, assume that $M$ is a non compact manifold and that $D$ is a bounded open subset of $M$ such that: i) the metric on $M\sm D$ is K\"ahler and $d$(bounded), and ii) $M\sm D$ has a strongly pseudo-convex boundary. Then the conclusion of Theorem \ref{kaehler} holds for this $M$.
\ee

{\small \textsc{Universit\'e de Nantes, D\'epartement de Math\'ematiques, 2 rue de la Houssini\`ere, BP 92208, 44322 Nantes cedex 03, France\\}}
e-mail: nader.yeganefar@math.univ-nantes.fr
\end{document}